\documentclass[graybox]{svmult}

\usepackage{type1cm}        
\usepackage{makeidx}         
\usepackage{graphicx}        
\usepackage{multicol}        
\usepackage[bottom]{footmisc}

\usepackage{newtxtext}       %
\usepackage{newtxmath}       

\usepackage{amsmath,mathtools,dsfont}
\usepackage[colorlinks=true, citecolor=red]{hyperref}
\newcommand{\Rn}{\mathbb{R}^n}
\newcommand{\R}{\mathbb{R}}
\newcommand{\N}{\mathbb{N}}

\newcommand{\Ec}{\mathcal{E}_s}
\newcommand{\Hc}{\mathcal{H}_s}
\newcommand{\Pc}{\mathcal{P}}

\makeindex             

\begin{document}

\title*{Long-range phase coexistence models with degenerate potentials.}


\author{Francesco~De~Pas, Serena~Dipierro, Enrico~Valdinoci}
\institute{Francesco~De~Pas \at University of Western Australia, 35 Stirling Highway, WA 6009 Crawley, Australia, \email{francesco.depas@uwa.edu.au}
\and Serena~Dipierro \at University of Western Australia, 35 Stirling Highway, WA 6009 Crawley, Australia, \email{serena.dipierro@uwa.edu.au}
\and Enrico~Valdinoci \at University of Western Australia, 35 Stirling Highway, WA 6009 Crawley, Australia, \email{enrico.valdinoci@uwa.edu.au}}


%
%
\maketitle

\abstract{
This survey offers an overview of recent advances in nonlocal phase transition problems, modeled by Ginzburg--Landau type energies of the form
\[
\frac{1}{4}\iint_{\R^{2n}\setminus (\R^n \setminus \Omega)^2}
\frac{|u(x)-u(y)|^2}{|x-y|^{n+2s}}\,dx\,dy
\;+\;
\int_\Omega W(u(x))\,dx.
\]
Here,~$W$ is a smooth and possibly \textit{degenerate} double well potential, with a polynomial control on its second derivatives near the wells.
The emphasis is on qualitative properties of minimizers and critical points of the energy functional.
}

\section{INTRODUCTION}\label{intro}
In this work we aim to shed light on several classical and recent results concerning an energy functional arising in phase transition models with long-range particle interactions. Specifically, given~$s\in(0,1)$ and a bounded set~$\Omega$ in~$\R^n$, we consider energies of the form
\begin{equation}\label{main_fun}
  \Ec(u;\Omega) := \Hc(u,\Omega)+\Pc(u,\Omega),
\end{equation}
where the nonlocal interaction term~$\Hc$ and the potential term~$\Pc$ are given, respectively, by
\begin{equation*}
\Hc(u,\Omega):=  \frac{1}{4}\iint_{\R^{2n}\setminus (\Rn \setminus \Omega)^2}\frac{|u(x)-u(y)|^2}{|x-y|^{n+2s}} \, dx\, dy,
\end{equation*}
and
\begin{equation*}\label{pot_term}
\Pc(u,\Omega) := \int_\Omega W(u(x)) \,dx.
\end{equation*}

Throughout this work, we assume that~$W$ is a smooth double-well potential with minima say at~$\pm1$, namely\footnote{As customary, when writing~$C^\vartheta(\Omega)$, if~$\vartheta>1$ we mean~$C^{k,\theta}(\Omega)$ with~$k\in\mathbb{N}$, $\theta\in(0,1]$, and~$\vartheta=k+\theta$.}
\begin{align}
&W \in C^{2,\beta}([-1,1]) \quad\mbox{for some}\quad \beta>0, \quad W>0 \ \mbox{in}\ (-1,1), \label{ojhuf}\\
&\mbox{and} \quad W(\pm1)=W'(\pm1)=0. \notag
\end{align}

We will also use the following notation, which is standard in the literature.
\begin{definition}
Let~$W$ be a potential satisfying~\eqref{ojhuf}. Then, we say that~$W$ is \emph{degenerate} if~$W''(\pm1)=0$ and~\emph{nondegenerate} if~$W''(\pm1)>0$.
\end{definition}

Functionals of the form~\eqref{main_fun} are known in literature as \textit{Ginzburg--Landau} type energies, in which the kinetic term~$\Hc$ is expressed through a fractional Sobolev seminorm.  
From a physical viewpoint, this implies that the interactions are not restricted to pairs of points both lying in~$\Omega$: points inside the domain interact directly with points in its complement.

Energies of this type arise naturally in several contexts. They provide effective descriptions of phase transition phenomena driven by nonlocal surface tension mechanisms (see e.g.~\cite{CozziValdNONLINEARITY, SV12, SV14, PSV13, CS14, MR3280032, MR4581189}) and play a central role in the Peierls--Nabarro theory for crystal dislocations (see e.g.~\cite{MR371203, MR1442163, MR2461827, MR2946964, GM12, DFV14, DPV15, MR3338445, BV16, MR3511786, MR3703556, MR4612096, MR4531940}).

As customary in many variational problems, particular attention has been devoted to equilibrium configurations of~\eqref{main_fun}. Such configurations correspond to solutions of the associated Euler--Lagrange equation, that goes in literature under the name of \emph{fractional Allen--Cahn equation}:
\begin{equation}\label{EL_eq}
L_s u(x) = W'(u(x)),
\end{equation}
where~$L_s$ denotes the fractional Laplacian of order~$s\in(0,1)$
\begin{equation*}\label{fralapc}
L_s u(x) := \text{\rm PV}_x \int_{\R^n} \frac{u(y)-u(x)}{|x-y|^{n+2s}}\, dy,
\end{equation*}
and~${\rm PV}$ stands for the Cauchy principal value.

In particular, among all equilibria, local minimizers of~\eqref{main_fun} play a distinguished role due to their enhanced stability properties. Nontrivial minimizers typically interpolate between the two pure phases of~$W$ and exhibit a rich qualitative structure, including symmetry and monotonicity properties, rigidity phenomena, and deep connections with nonlocal minimal surfaces (see e.g.~\cite{MR2177165, MR2498561, MR2644786, SV12, MR3035063, MR3148114, SV14, BV16, MR3740395, MR3812860, MR3939768, MR4124116, MR4050103, MR4938046}). For these reasons, their analysis constitutes a central topic in the theory of (classical and) nonlocal partial differential equations.

In the literature, such minimizers are commonly referred to as \emph{layer solutions} or \emph{transition layers}, a terminology that we shall adopt throughout this work. More precisely, we introduce the following notion.

\begin{definition}
Let us consider the class of admissible functions
\begin{equation*}
\mathcal{X} := \left\{ f \in L^1_{\rm loc}(\R) \ \text{s.t.} \ \lim_{x\to\pm\infty} f(x)=\pm1 \right\},
\end{equation*}
and let~$u\in\mathcal{X}$ be a solution to~\eqref{EL_eq}. We say that~$u$ is a \textit{layer solution} if it is monotone increasing from~$-1$ to~$+1$ in one of the spatial variables.
\end{definition}

The survey is organized as follows. In Section~\ref{history}, we review the main classical results on the topic, providing a broad perspective on related questions and possible connections to other problems in phase transitions and nonlocal variational models. Section~\ref{p0op} is devoted to the analysis of monotone heteroclinic connections in dimension $n=1$, highlighting both existence and qualitative properties of these solutions. Finally, Section~\ref{8f5} presents an overview of the reconstruction of a potential $W$ from a prescribed boundary layer profile in dimension $n=1$.

\section{Classical results and background}\label{history}
In this paragraph we provide a brief overview of some of the main contributions to the Allen--Cahn equation developed over the years. As mentioned above, the study of qualitative properties of monotone minimizers for energies of this type has a long and well-established history in the literature.

In particular, this line of research traces back to a conjecture formulated by De Giorgi in the late~$1970$s, originally stated in the local setting. The conjecture can be formulated as follows.

\begin{problem}
Let~$n\leq 8$ and~$u$ be a bounded entire solution of
\begin{equation}\label{epi90}
-\Delta u = W'(u),
\end{equation}
with~$\partial_{x_1}u>0$ in~$\R^n$. 

Is it true that~$u$ is one-dimensional? Namely, does there exist~$e\in S^{n-1}$ and a function~$u_0:\R\to \R$ such that~$u(x) = u_0(e\cdot x)$ for any~$x\in\R^n$?
\end{problem}

Equivalently, the conjecture asks whether the level sets of such monotone solutions must necessarily be hyperplanes. De Giorgi’s conjecture was proved to be true for~$n\leq 3$ see \cite{MR1655510,MR1775735}. The cases~$4\leq n\leq 8$ were settled under the additional assumption that
\begin{equation*}
\lim_{x_n\to\pm\infty}u(x',x_n)=\pm1 \quad\mbox{for any } x'\in\R^{n-1},
\end{equation*}
as shown in~\cite{MR2480601}. For~$n\geq 9$, a counterexample was constructed in~\cite{MR2473304}, showing that the dimensional threshold in the conjecture is sharp.

Analogous symmetry questions have also been investigated in the nonlocal setting, where the Laplacian in~\eqref{epi90} is replaced by the fractional operator~$L_s$. In this framework, the conjecture has been proved in dimension $n=2$ in~\cite{MR2498561}. 
In dimension $n=3$, it was established for $s\geq \tfrac12$ in~\cite{MR2644786, MR3148114}, 
and later extended to the full range $s\in(0,1)$ in~\cite{MR4124116, MR3740395}. 
The case $n=4$ and $s=\tfrac12$ was treated in~\cite{MR4050103}. 
Moreover, the validity of the conjecture for a broad class of truncated kernels was proved in~\cite{MR3610941}.

Another key aspect of the Allen--Cahn equation is its deep connection with minimal surface theory. This relationship was first unveiled in the seminal work of Modica and Mortola in~\cite{MR445362} for the local case, where the authors identified the sharp--interface limit of the Ginzburg--Landau energy. Later, in~\cite{MR1310848, farina}, uniform density estimates for minimizers of the Allen--Cahn energy were established.

These ideas were subsequently extended to the nonlocal setting. In particular, in~\cite{SV12} the authors investigated the connection between solutions of the fractional Allen--Cahn equation and (non)local minimal surfaces, obtaining a fractional analogue of the classical $\Gamma$-convergence result of Modica and Mortola. 
Subsequently, uniform density estimates for solutions and minimizers in the fractional framework were proved in~\cite{SV14}, playing a role analogous to their local counterparts in rigidity and classification results.

For completeness, we also refer the reader to~\cite{AB94, PV20} for $\Gamma$-convergence results in one dimension, corresponding respectively to the cases $s=\tfrac12$ and $s\in(0,1)$. 

More recently, density estimates in the fractional setting have been established for~\emph{degenerate} double-well potentials in~\cite{DFGV26, DFGV25}. 

Altogether, these results show that the geometry of the level sets of solutions to the Allen--Cahn equation reflects, at the interface, the structure of minimal surfaces . From this perspective, rigidity results for level sets—such as those predicted by De Giorgi’s conjecture for monotone solutions—can be interpreted as diffuse counterparts of classical rigidity phenomena in minimal surface theory.

\section{Monotone heteroclinic connections in dimension $n=1$}\label{p0op}

As discussed in Sections~\ref{intro} and~\ref{history}, transition layers connecting stable states of a double-well potential exhibit a rich and well-structured behavior. In this section, we collect a number of results concerning existence, regularity, and decay properties of such minimizers in the one-dimensional setting $n=1$.

We begin by clarifying the notion of minimizers for the energy functional~\eqref{main_fun}.

\begin{definition}\label{defini}
Let~$\Omega$ be a bounded domain of~$\R$. A measurable function~$u : \R \to \R$ is called a {\rm local minimizer} of~$\Ec$ in~$\Omega$ if~$\Ec(u;\Omega)<+\infty$ and
\[
\Ec(u;\Omega) \leqslant \Ec(u+\phi;\Omega),
\]
for every~$\phi\in C^\infty_0(\Omega)$.

Moreover, we say that a measurable function~$u : \R \to \R$ is a {\rm class~A minimizer} of~$\Ec$ if it is a local minimizer in every bounded domain~$\Omega\subset \R$.
\end{definition}

The introduction of class~A minimizers is dictated by a structural feature of the problem itself. Indeed, due to the nonlocal nature of the energy~$\Ec$, its value may diverge when computed over unbounded domains. For this reason, it is customary to keep track of the growth rate of~$\Ec$ via the auxiliary functions
\begin{equation*}\label{fnc_psi}
    \Psi_s(\rho):=
    \begin{cases}
        \rho^{1-2s}, & \mbox{if } s \in (0,1/2),\\
        \log \rho, & \mbox{if } s=1/2,\\
        1, & \mbox{if } s \in (1/2,1).
    \end{cases}
\end{equation*}
and 
\begin{equation*}
\mathcal{G}(u):= \limsup_{\rho \to +\infty}\frac{\Ec(u;[-\rho,\rho])}{\Psi_s(\rho)}.
\end{equation*}

This scaling naturally arises from~\cite[Proposition~$3.1$]{CozziValdNONLINEARITY}, where it is shown that any local minimizer~$u\in\mathcal{X}$ of the energy~$\Ec(\cdot,B_{R+2})$ satisfies
\[
\Ec(u,B_R) \leq C \Psi_s(R),
\]
for any~$R\geq 3$.

A first fundamental result on one-dimensional heteroclinic connections is due to~\cite{PSV13}. The authors consider a~\emph{nondegenerate} double-well potential satisfying~\eqref{ojhuf} and establish  the following estimates for boundary layers.

\begin{theorem}[Theorem~$2$ in~\cite{PSV13}]\label{cgd8}
Let~$n=1$. There exists a unique (up to translations) nontrivial class~A minimizer~$u \in \mathcal{X}$ of the energy~$\Ec$, which is strictly increasing. Moreover,~$u$ solves~\eqref{EL_eq} and is unique (up to translations) also within the class of monotone solutions of~\eqref{EL_eq}.

In addition,~$u\in C^2(\R)$ and there exists a constant~$C\geq 1$ such that
\[
|u(x)-{\rm sign}(x)|\leq C |x|^{-2s}
\quad\mbox{and}\quad
|u'(x)|\leq C |x|^{-1-2s},
\]
for all sufficiently large~$|x|$. Finally,~$\mathcal{G}(u)<+\infty$.
\end{theorem}
See also~\cite{MR3280032} for related results.

A further step in the theory is achieved in~\cite{CP16}, where the authors consider a more general class of nonlocal energies~$\Ec$. In particular, the Sobolev seminorm in~$\Hc$ is replaced by
\begin{equation}\label{fu7}
\mathcal{H}_K (u,\Omega) := \frac14 \iint_{\R^{2}\setminus(\R\setminus\Omega)^2} |u(x)-u(y)|^2 K(x-y)\, dx\, dy,
\end{equation}
where~$K$ is a measurable, positive and symmetric kernel modeled on that of the fractional Laplacian. More precisely, the kernel~$K$ is assumed to satisfy
\[
\exists\, 0<\lambda\leq\Lambda,\ r_0>0:\quad
\frac{\lambda\,\mathds{1}_{B_{r_0}}(z)}{|z|^{n+2s}} \leq K(z) \leq \frac{\Lambda}{|z|^{n+2s}}.\footnote{
This hypothesis is very broad, it allows for instance kernels of the form
\[
K(z):=\mathds{1}_{B_{r_0}}(z)\frac{a(z)}{|z|^{n+2s}},
\]
with~$a$ bounded and strictly positive .
}
\]

In addition, in~\cite{CP16} the authors work with a symmetric potential~$W$, which yields symmetric minimizers. As a consequence,~\cite[Theorem~$2$]{CP16} provides an analogue of Theorem~\ref{cgd8} within this more general framework.

The most recent developments in this direction aim at substantially relaxing the assumptions on the potential $W$, in particular allowing for a \emph{degenerate} behavior at the wells~$\pm1$. 
By permitting the potential to approach its wells at a polynomial rate, one can treat a broader class of potentials, and it is natural to expect that changes in the behavior of $W$ (and its derivatives) influence the structure and properties of the transition layers.
At the same time, full generality is preserved in the choice of the kernel defining the elastic contribution $\mathcal{H}_K$ as in~\eqref{fu7}.
 Such results can be found in~\cite{OURREC, OURREC2}. In particular, the potential~$W:\R\to[0,+\infty)$ is assumed to satisfy the conditions in~\eqref{ojhuf} and that there exist~$C_2 \geq C_1 >0$, $C_4\ge C_3>0$, $\xi \in (0,1)$, $\alpha \geq \beta \geq 2$, and~$\gamma \geq \delta \geq 2$ such that
\begin{equation*}\label{pot_deg}
\begin{cases}
C_1 (1+t)^{\alpha-2} \leq W''(t) \leq C_2 (1+t)^{\beta-2}
& \text{for } t \in (-1,-1+\xi],\\
C_3 (1-t)^{\gamma-2} \leq W''(t) \leq C_4 (1-t)^{\delta-2}
& \text{for } t \in [1-\xi,1).
\end{cases}
\end{equation*}

We emphasize that condition~\eqref{pot_deg} is very general: it allows the potential~$W$ to be~\emph{degenerate} and even to exhibit oscillatory behavior near the wells. Moreover, the symmetry assumption imposed in~\cite{CP16} is removed.

Under these assumptions, new and sharp decay estimates for one-dimensional minimizers are obtained.

\begin{theorem}[Theorem~$1.4$ in~\cite{OURREC}]\label{bou}
Let~$n=1$ and assume that
\begin{equation}\label{riar}
\max \{ \alpha-\beta, \gamma-\delta \} <1.
\end{equation}
Then, within the class~$\mathcal{X}$, there exists a unique (up to translations) nontrivial class~A minimizer~$u$ of~$\Ec$, satisfying~$\mathcal{G}(u)<+\infty$. 

Moreover,~$u$ is strictly increasing, belongs to~$C^{1+2s+\theta}(\R)$ for some~$\theta\in(0,1)$, and is the unique increasing solution of
\[
L_K u = W'(u) \quad \text{in } \R.
\]

In addition, there exist constants~$C_1$, $C_2>0$ and~$R>0$ such that the following asymptotic estimates hold:

\begin{equation}\label{asymp_decay}
\left\{
\begin{alignedat}{2}
1+u(x) &\le C_2 |x|^{-\frac{2s}{\alpha-1}} &\qquad& \text{for } x \le -R,\\
1-u(x) &\le C_2 |x|^{-\frac{2s}{\gamma-1}} &      & \text{for } x \ge R.
\end{alignedat}
\right.
\end{equation}

\begin{equation}\label{398r7gree}
\left\{
\begin{alignedat}{2}
u'(x) &\ge C_1 |x|^{-\left(1+\frac{2s(\alpha-\beta+1)}{\alpha-1}\right)}
&\qquad& \text{for } x \le -R,\\
u'(x) &\ge C_1 |x|^{-\left(1+\frac{2s(\gamma-\delta+1)}{\gamma-1}\right)}
&      & \text{for } x \ge R.
\end{alignedat}
\right.
\end{equation}

\begin{equation}\label{asymp_decay_lowbound}
\left\{
\begin{alignedat}{2}
1+u(x) &\ge C_1 |x|^{-\frac{2s}{\beta-1}} &\qquad& \text{for } x \le -R,\\
1-u(x) &\ge C_1 |x|^{-\frac{2s}{\delta-1}} &      & \text{for } x \ge R.
\end{alignedat}
\right.
\end{equation}

\begin{equation}\label{eq:asymp-derivata}
\left\{
\begin{alignedat}{2}
u'(x) &\le C_2 |x|^{-\left(1+\frac{2s(\beta-\alpha+1)}{\beta-1}\right)}
&\qquad& \text{for } x \le -R,\\
u'(x) &\le C_2 |x|^{-\left(1+\frac{2s(\delta-\gamma+1)}{\delta-1}\right)}
&      & \text{for } x \ge R.
\end{alignedat}
\right.
\end{equation}
\end{theorem}

We mention that hypotheses~\eqref{riar} can be interpreted as a small-amplitude oscillation condition for~$W$ near the pure phases. 
Finally, combining~\cite[Remark~$1.7$]{DPDV} and~\cite[Remark~$1.6$]{OURREC}, one obtains the following optimality result.

\begin{theorem}
If~$\alpha=\beta$ and~$\gamma=\delta$, then the estimates~\eqref{asymp_decay}, \eqref{398r7gree}, \eqref{asymp_decay_lowbound}, and~\eqref{eq:asymp-derivata} are optimal. If instead the potential is oscillatory (that is,~$\alpha\neq\beta$ or~$\gamma\neq\delta$), then~\eqref{asymp_decay}, \eqref{398r7gree}, and~\eqref{asymp_decay_lowbound} are optimal.
\end{theorem}

\begin{svgraybox} \textbf{Further investigation}:\\
It would be interesting to investigate whether the upper bounds in~\eqref{eq:asymp-derivata} remain optimal even in the oscillatory case~$\alpha\neq\beta$ and~$\gamma\neq\delta$.
\end{svgraybox}

\section{Reconstructing a potential from the layer in~$n=1$}\label{8f5}
As recalled in Section~\ref{p0op}, the classical approach to phase transition problems consists in fixing a double-well potential~$W$, typically chosen on phenomenological or heuristic grounds, and then studying the existence, regularity, and qualitative properties of the associated minimizers. In this perspective, results such as Theorems~\ref{cgd8} and~\ref{bou} show that sufficiently regular double-well potentials give rise to transition layers solving the Allen--Cahn equation~\eqref{EL_eq}, and that the behavior of~$W$ and its derivatives near the wells governs the asymptotic decay of the layer.

This naturally raises the converse question, which lies at the core of this subsection.

\begin{problem}\label{xyd6}
Let~$n=1$. Given a transition layer whose asymptotic behavior is prescribed a priori, is it possible to construct a double-well potential~$W$ such that~\eqref{EL_eq} is satisfied? If so, how does the shape of the layer determine the behavior of~$W$ and its derivatives near the wells?
\end{problem}

This inverse viewpoint is consistent with the phenomenological nature of phase transition models. In many applications, the potential is not derived from first principles but postulated in the spirit of Landau's theory for phase transition models, where the free energy is expanded in powers of an order parameter. From this perspective, reconstructing the potential from observable features of an interface—such as decay rates toward the pure phases—appears both natural and conceptually meaningful.

The relevance of Problem~\ref{xyd6} is already suggested by the works~\cite{CS14, MR3280032}. Indeed, their results show that the mere existence of a monotone transition layer imposes strong structural constraints on the potential.

\begin{theorem}[Theorem~$2.4$ in~\cite{MR3280032}]\label{jbvr9}
Let~$f\in C^{1,\gamma}(\R)$ with~$\gamma>\max\{0,1-2s\}$ and let~$W'=f$. Then there exists a monotone increasing solution~$u\in\mathcal{X}$ of~\eqref{EL_eq} if and only if
\[
W'(\pm1)=0 \quad\mbox{and}\quad W>W(1)=W(-1) \ \mbox{in } (-1,1).
\]
If, in addition,~$W''(-1)>0$ and~$W''(1)>0$, then such a solution is unique up to translations.
\end{theorem}

This theorem shows that the existence of a layer solution forces the potential to be of double-well type, thus providing a first partial answer to Problem~\ref{xyd6}.

An alternative perspective arises in the study of crystal dislocation models. 
In~\cite{DPV15}, the authors explicitly build a potential starting from a prescribed transition layer. More precisely, they define an even and positive function~$w\in C^\infty(\R)$ by
\[
w(x):=\frac1{|x|^{1+2s}}\qquad\mbox{for }|x|\geq1.
\]
Then, they let~$A:=\|w\|_{L^1(0,+\infty)}$ and introduce the function~$\phi\in C^\infty(\R)$ defined as
\[
\phi(x):=\frac1A\int_0^x w(t)\,dt.
\]
They then set
\begin{equation*}\label{i3456789097654bvcxs}
\begin{split}
& g(t):=L_s\phi(t)\quad\text{for all }t\in\R,\\
& h(r):=g(\phi^{-1}(r))\quad\text{for all }r\in(-1,1),
\end{split}
\end{equation*}
and define the potential
\begin{equation}\label{pmes}
V(r):=\int_{-1}^r h(\rho)\,d\rho \quad\mbox{for any }r \in (-1,1).
\end{equation}

The resulting potential enjoys the following properties.

\begin{theorem}[Propositions~$6.2$ and~$6.3$ in~\cite{DPV15}]\label{cut8}
The function~$\phi$ and the potential~$V$ satisfy
\[
L_s\phi = V'(\phi)\quad\mbox{in }\R.
\]
Moreover,~$V\in C^\infty(-1,1)$,~$V$ is symmetric with respect to~$0$. In addition
\[
V(\pm1)=V'(\pm1)= 0,\quad V(r)>0 \quad\mbox{for any} \ r\in(-1,1)  \quad\mbox{and}\quad V''(\pm1)=2A.
\]
\end{theorem}

Beyond its intrinsic interest, this construction has concrete analytical applications. In particular, the function~$\phi$ is used as a barrier to refine the decay estimates for solutions to~\eqref{EL_eq} in the case of~\emph{ nondegenerate} double-well potentials. In particular, Theorem~\ref{cut8} leads to

\begin{theorem}[Proposition~$7.2$ in~\cite{DPV15}]\label{KVV}
Let~$n=1$ and~$u$ solve~\eqref{EL_eq} with~$s\in(1/2,1)$ and~$W$ be nondegenerate. Then there exist constants~$a_0>0$ and~$C>0$ such that
\[
\left|u(x)-H(x)+\frac{a_0}{2s}\frac{x}{|x|^{1+2s}}\right|\leq\frac{C}{|x|^{1+2s}}.
\]
\end{theorem}

\begin{remark}
The case~$s\in(0,1/2)$ corresponding to Theorem~\ref{KVV} has been addressed in~\cite[Theorem~$1.1$]{DFV14}.
\end{remark}

Despite its effectiveness, this example leaves several aspects of Problem~\ref{xyd6} open. First, the construction is restricted to symmetric layers, a setting in which stronger conclusions are typically available. It is therefore natural to ask whether a non-symmetric layer~$\phi$ would still generate a (possibly non-symmetric) double-well potential. Second, Theorem~\ref{cut8} does not provide detailed information on the regularity of~$V$ at the wells: in principle,~$V''$ could not even be H\"older continuous near the wells.

A more comprehensive answer to Problem~\ref{xyd6} is provided in~\cite{OURREC2}, where this topic is treated as the main object of study. The authors there consider polynomial-type layers designed as follows. Let~$\alpha,\beta\in(0,2s]$,~$\kappa>0$, and~$C_1,C_2>0$, and assume that~$\phi\in C^\infty(\R)$ satisfies~$\phi'>0$ and
\begin{equation}\label{valesem}
\phi(x):=
\begin{cases}
-1+C_1|x|^{-\alpha} & \text{if } x<-\kappa,\\
1-C_2|x|^{-\beta} & \text{if } x>\kappa.
\end{cases}
\end{equation}
The potential~$V$ is then defined as in~\eqref{pmes}.

\begin{theorem}[Theorem~$2.1$ in~\cite{OURREC2}]\label{minth}
The function~$\phi$ and the potential~$V$ satisfy
\[
L_s\phi = V'(\phi)\quad\mbox{in }\R.
\]
Moreover,~$V\in C^\infty(-1,1)$ and
\begin{equation*}\label{48397tasvjkdfbgkewguo}
V(\pm1)=0 \quad\mbox{and}\quad V(r)>0 \quad\mbox{for all } r\in(-1,1).
\end{equation*}
In addition,
\begin{equation*}\label{ham}
\lim_{r\to-1^+}\frac{V(r)}{(1+r)^{\frac{2s}\alpha+1}}
=\frac{\alpha C_1^{-\frac{2s}{\alpha}}}{(2s+\alpha)s},
\qquad
\lim_{r\to1^-}\frac{V(r)}{(1-r)^{\frac{2s}\beta+1}}
=\frac{\beta C_2^{-\frac{2s}\beta}}{(2s+\beta)s},
\end{equation*}
and
\begin{equation*}\label{hamm}
\lim_{r\to-1^+}\frac{V'(r)}{(1+r)^{\frac{2s}\alpha}}
=\frac{C_1^{-\frac{2s}{\alpha}}}{s},
\qquad
\lim_{r\to1^-}\frac{V'(r)}{(1-r)^{\frac{2s}\beta}}
=-\frac{C_2^{-\frac{2s}\beta}}{s}.
\end{equation*}

Furthermore, for any~$i\in\N$,
\begin{equation*}\label{0ijuhed}
\lim_{r\to-1^+}\frac{V^{(i+1)}(r)}{(1+r)^{\frac{2s}\alpha-i}}\in(-\infty,+\infty),
\qquad
\lim_{r\to1^-}\frac{V^{(i+1)}(r)}{(1-r)^{\frac{2s}\beta-i}}\in(-\infty,+\infty).
\end{equation*}

If~$i\in\N\setminus\{0\}$ and~$2s\geq\alpha i$, then
\begin{equation}\label{trep1}
\lim_{r\to-1^+}\frac{V^{(i+1)}(r)}{(1+r)^{\frac{2s}\alpha-i}}
=\frac1s C_1^{-\frac{2s}{\alpha}}
\prod_{j=0}^{i-1}\left(\frac{2s}\alpha-j\right).
\end{equation}
Similarly, if~$i\in\N\setminus\{0\}$ and~$2s\geq\beta i$, then
\begin{equation}\label{trep2}
\lim_{r\to1^-}\frac{V^{(i+1)}(r)}{(1-r)^{\frac{2s}\beta-i}}
=\frac{(-1)^{i+1}}{s} C_2^{-\frac{2s}\beta}
\prod_{j=0}^{i-1}\left(\frac{2s}\beta-j\right).
\end{equation}

Finally\footnote{
As usual, we use the notation~$\lfloor x\rfloor:=\max\{y\in\N:\ y\le x\}$.
},
\begin{equation*}\label{REgg}
V\in C^{\lfloor\frac{2s}\alpha\rfloor,1}([-1,0])
\quad\mbox{and}\quad
V\in C^{\lfloor\frac{2s}\beta\rfloor,1}([0,1]).
\end{equation*}
\end{theorem}

Theorem~\ref{minth} thus provides a complete answer to Problem~\ref{xyd6} for polynomial-type layers of the form~\eqref{valesem}. It establishes a precise correspondence between the decay rates of~$\phi$ at infinity and the order of vanishing of~$V$ and its derivatives at the wells.

In particular, choosing a polynomially decaying layer produces a potential that exhibits a power-type structure near~$\pm1$. We stress that the fact that the fractional Laplacian of the transition layer enjoys good scaling properties on derivatives is nontrivial given the 
nonlocal nature of the operator\footnote{In spite of its own scale invariance,
the nonlocal character of the fractional Laplacian operator does not make it compatible, in general, with the notion of power-like functions. To see this, given~$a<b$ and~$\lambda\in\R$,
we observe that, on the one hand, if~$f(x)=x^\lambda$ for all~$x\in(a,b)$, then all the derivatives of~$f$ are power-like in~$(a,b)$. On the other hand, 
by~\cite[Theorem~1.5]{MR4297378},
for every~$m\in\N$, every~$F\in C^m(\R)$, and
every~$\epsilon>0$, one can find functions~$f_\epsilon$ and~$\eta_\epsilon$ such that~$\|f-f_\epsilon\|_{C^m((a,b))}\le\epsilon$,
$\|\eta_\epsilon\|_{L^\infty((a,b))}\le\epsilon$, and
$$(-\Delta)^s f_\epsilon=F(x)+\eta_\epsilon\qquad{\mbox{in }}(a,b).$$
For example, the fractional Laplacian of a function which ``looks like a power in $(a,b)$'',
may well ``look like an exponential in~$(a,b)$''.}. 
Each differentiation of~$V$ lowers the effective order by one, as made explicit by~\eqref{trep1} and~\eqref{trep2}. The case~$\alpha=\beta=2s$ corresponds to a~\emph{nondegenerate} potential, while~$\alpha,\beta<2s$ yields \emph{degenerate} wells, in agreement with the framework of~\cite{DPDV}. More generally, for any~$i\in\N$, selecting~$\alpha,\beta\in(0,2s/i]$ guarantees that~$V\in C^{i,1}((-1,1))$ and~$V^{(i)}(\pm1)=0$.

\begin{svgraybox} \textbf{Further investigation}:\\
It would be interesting to approach Problem~\ref{xyd6} for transition layers exhibiting asymptotic behaviors different from the polynomial ones considered above. From a modeling perspective, non-polynomial decays may encode distinct physical regimes and lead to different behaviors of~$W$ the wells, even when such decays are asymptotically comparable to powers of~$x$. 
A paradigmatic example is provided by the function 
\[
u(x):=\frac2\pi \arctan(x) \quad\mbox{and}\quad V(\rho):=\frac1{\pi^2}\big(\cos(\pi\rho)-\cos(\pi)\big).\]
It is known that (see~\cite{OURREC2} and~\cite[Appendix~L]{AV19} )
\[
L_{1/2}u(x)=V'(u(x))\quad\mbox{in }\R.
\]
This example shows that the inverse construction problem admits solutions beyond the polynomial setting. Extending the analysis of Problem~\ref{xyd6} appears therefore potentially rich in new phenomena.
\end{svgraybox}

\begin{svgraybox} \textbf{Further investigation}:\\
It would be interesting to exploit Theorem~\ref{minth} to construct suitable barriers for solutions of~\eqref{EL_eq} in the case of~\emph{degenerate} or oscillatory potentials, with the aim of obtaining an analogue of Theorem~\ref{KVV} in the~\emph{degenerate} setting.
\end{svgraybox}

%
%
%

\end{document}